\documentclass[12pt]{amsart}

\usepackage{fancyhdr}
\usepackage{algorithmic}
\usepackage{algorithm}

\theoremstyle{plain}
\newtheorem{lemma}{Lemma}[section]
\newtheorem{teo}[lemma]{Theorem}
\newtheorem{propo}[lemma]{Proposition}

\theoremstyle{definition}

\newtheorem{example}[lemma]{Example}

\theoremstyle{remark}

\newcommand{\z}{\mathbb{Z}}

\newcommand{\rt}{\rightarrow}
\newcommand{\e}{\epsilon}

\begin{document}
\setlength{\baselineskip}{15pt}

\title{Representations of matroids}
\author{Massimiliano Lunelli and Antonio Laface}

\begin{abstract}
In this paper we give a necessary and sufficient criterion for  
representability of a matroid over an algebraic closed field. 
This leads to an algorithm, based on an extension of Gr\"{o}bner Bases, 
in order to decide if a given matroid is representable over such a field.
\end{abstract}
\keywords{matroids, Grobner basis, representability}
\maketitle

\section{Introduction}
A matroid $M$ on a finite set $S$ is said to be representable over a field $k$
if there exists a vector space $V$ over $k$ and an injection $\phi\colon
S\rightarrow V$ such that a subset $X$ of $S$ is independent in $M$ if and
only if the vectors of $\phi(X)$ are linearly independent over $k$. \\
The problem of the existence of such representations has been largely studied
or fields of characteristic $2,3$ and $4$. In these cases, the main result is
that there exists a finite list of non-representable matroids, such that, for
each matroid $M$ having a minor in the list, $M$ is not representable. \\
There are also results, connecting representability in finite characteristic,
with representability in characteristic $0$.\\
In this paper the point of view is quite different; an algorithm, based on
Gr\"{o}bner bases, is given, in order to establish if there exists a linear
representation for a given matroid over some field. 
This transform the problem of representability over such fields in a purely 
algebraic question. \\
There are two problems in this kind of approach: the first one is computational,
due to the large amount of calculus involved in the determination of the
Gr\"{o}bner bases. This step may be exceeded by finding an appropriate algebraic
presentations for the given matroid. This question leads to a sufficient algorithm for the representability of a matroid.
The second obstacle is theoric, it depends from the fact that a-priori it is 
not known over which field a representation may be given. 
For example, there are matroids that are representable only over a field 
of characteristic $2$, for such a matroid, the response of a representability
test will be negative if the field has a characteristic different from $2$.
This lead us to develop a Gr\"{o}bner bases algorithm over $\z$.

\section{Representability and Algorithms}

In order to find a representation of the matroid $M$, we are lead to construct
an $r\times n$ matrix, whose entries belong to some unknown field $k$:

\begin{equation}
\label{mat1}
\left(
\begin{array}{ccc}
x_{1,1} & \cdots & x_{1,n} \\
\vdots & & \vdots \\
x_{r,1} & \cdots & x_{r,n} \\
\end{array}
\right)
\end{equation}

Consider the set of bases of $M$. If $v_{\sigma(1)},\cdots,v_{\sigma(r)}$ is one
such a basis ($\sigma$ is a choose of $r$ numbers in $\{1\cdots n\}$), then take
the determinant of the corresponding $r\times r$ minor. \\
Let $\bar{x}=x_{1,1}\cdots x_{r,n}$ and let $P_1(\bar{x})\cdots P_s(\bar{x})$ 
be the polynomials obtained from all the bases of $M$.
Let $Q_1(\bar{x})\cdots Q_l(\bar{x})$ be polynomials obtained by
the determinant of the minors arising from the circuits and let $I$ be the ideal
of $k[x_{1,1}\cdots x_{r,n}]$ generated by the $Q_i$. 
With this notation, we have the following:

\begin{teo}
\label{main}
A matroid $M$ of cardinality $n$ and rank $r$ is representable over an algebraic
closed field $k$ if and only if
\[
\prod_{i=1}^sP_i(\bar{x})\not\in \rm{Rad}(I)
\]
where $<Q_1(\bar{x}),\cdots, Q_l(\bar{x})>$ is the ideal of $k[\bar{x}]$
generated by $Q_1,\cdots, Q_l$.
\end{teo}
\begin{proof}
To find a representation of $M$ over $k$ is equivalent to find  
$\bar{x}_0$ such that $Q_1(\bar{x}_0)=\cdots =Q_l(\bar{x}_0)=0$
and $P_i(\bar{x}_0)\neq 0$ for each $i\in\{1,\cdots, s\}$.
This means that the hypersurface given by the zero locus of
$\prod_{i=1}^sP_i(\bar{x})=0$ does not contain the variety given by the
intersection of the $Q_i$. By Hilbert's Nullstellensatz, this is equivalent to
ask that $\prod_{i=1}^sP_i(\bar{x})$ does not belong to the radical of the
ideal $I$ generated by the $Q_i$'s.
\end{proof}

In order to apply the preceding proposition we need to solve an ideal 
membership problem. This may be done by using Gr\"{o}bner basis. 
First of all, we may suppose that the first $r$ elements of $M$ are a basis of
the matroid, and hence in matrix \ref{mat1} we may take the first $r\times r$
minor to be the identity matrix. 

\begin{propo}
\label{pro:circ2}
A matroid $M$ is completely determined by the set of circuits of order less 
than or equal to its rank.
\end{propo}
\begin{proof}
Each basis is determined by this set, since if a set of $r$ elements 
($r=\rm{rank}(M)$) is dependent, then it must contain a circuit with $l\leq r$
elements.
\end{proof}

This allows us to restrict our attention to the $Q_i(\bar{x})$ coming from such
circuits. Equations coming from circuits with more than $\rm{rank}(M)$ 
elements are automatically satisfied when one search a representation in a
vector space of dimension $\leq\rm{rank}(M)$.

The first algorithm, allows one to determine if a matroid $M$ is representable 
over a given algebraically closed field $k$.

\begin{algorithm}
\caption{{\em Representation of matroids over an algebraically closed field}}
\label{alg:field-ns}

\begin{algorithmic}[1]

\REQUIRE $Q_1\cdots Q_l, P_1\cdots P_r$, the field $k$ and $\bar{x}=[\cdots x_{i,j}\cdots]$

\STATE $I := <Q_1\cdots Q_l>, P := \prod P_i, p={\rm char}(k)$
\STATE {\bf Find} a Gr\"{o}bner basis of $I$ over $\z_p[\bar{x}]$ 
\IF{$P \equiv 0\pmod{\rm{Rad}(I)}$}
\STATE The matroid in not representable over $k$.
\ENDIF

\end{algorithmic}
\end{algorithm}

Step 3 may be hard to finish, due to the great amount of calculation (the polynomial $P$ may be 
very bigger!). In order to avoid this obstacle it is possible to test if $P_i \equiv 0\pmod{I}$
for each bases equation $P_i$. But this idea does not work well, since in almost all cases it 
happens that some product of the $P_i$ belong to $\rm{Rad}(I)$ even if no one of the $P_i$ 
belong to it. The following proposition suggests another system.

\begin{propo}
\label{pro:circ1}
Given a basis for $M$, there exists only one circuit containing a given 
element and elements from the basis. 
\end{propo}
\begin{proof}
Clearly one such circuits must exists. Suppose that there is another circuit
with this property. Then the element $v$ may be expressed in two way as a linear
combination of elements of the given basis. By substitution, this would give a 
linear dependence between the basis elements.
\end{proof}

This means that it is possible to establish which of the $x_{i,j}$ vanish and which not: 
it is sufficient to consider the circuit given by the first $r$ columns and by the column
containing $x_{i,j}$. For each one of the non vanishing $x_{i,j}$, there is a corresponding 
bases B, such that the corresponding polynomial is $x_{i,j}$.
This means that if we ask that the $x_{i,j}$ do not vanish, then we already take in account 
many bases equation. The easiest way to do this is to change the ideal generated by the $Q_i$,
by adding another term: $1-t\prod x_{i,j}$.
This lead to the following:

\begin{algorithm}
\caption{{\em Representation of matroids over an algebraically closed field}}
\label{alg:field-s}

\begin{algorithmic}[1]

\REQUIRE $Q_1\cdots Q_l, P_1\cdots P_r$, ${\rm char}(k)$ and $\bar{x}=[\cdots x_{i,j}\cdots]$

\STATE $I := <Q_1\cdots Q_l,1-t\prod x_{i,j}>, p={\rm char}(k)$
\STATE {\bf Find} a Gr\"{o}bner basis of $I$ over $\z_p[\bar{x},t]$ 

\FORALL{$P\in\{P_1\cdots P_r\}$} 

 \IF{$P \equiv 0\pmod{I}$}
  \STATE The matroid in not representable over $k$.
 \ENDIF

\ENDFOR

\end{algorithmic}
\end{algorithm}

It is possible to consider the problem of the representation of the given 
matroid over every field. 
The idea is to consider Gr\"{o}bner bases of polynomials over the ring of 
integers. 
The following algorithm gives a sufficient condition for non-representability
of matroids over every field.

\begin{algorithm}
\caption{{\em Representation of matroids over any field}}
\label{alg:main}

\begin{algorithmic}[1]

\REQUIRE $Q_1\cdots Q_l, P_1\cdots P_r, L$ and $\bar{x}=[\cdots x_{i,j}\cdots]$

\STATE $I := <Q_1\cdots Q_l,1-t\prod x_{i,j}>, L := \emptyset$

\STATE {\bf Find} a Gr\"{o}bner basis $G$ of $I$ over $\z[\bar{x},t]$ 
     \FORALL{division by $n$ performed to obtain $G$}
     \FORALL{$p\ {\rm prime}, p\mid n$}
       \STATE use {\bf Algorithm 2} to test representability in ${\rm char}(p)$ 
     \ENDFOR
     \ENDFOR
        
  \IF{M is not representable over the preceding fields and $1\in I$}
    \STATE The matroid in not representable over any field
  \ENDIF

 \FORALL{$P\in\{P_1\cdots P_r\}$} 

  \IF{M is not representable over the preceding fields and $P \equiv 0\pmod{I}$} 
   \STATE The matroid in not representable over any field
  \ENDIF
 \ENDFOR

\end{algorithmic}
\end{algorithm}

\section{Examples}

As an example of non-representable matroid, we consider the non-Pappus matroid,
in this case there is a contradiction, between the conditions imposed 
by the circuits and one base.\\

\begin{example}
\label{ex:non-pappus}

Consider the non-Pappus matroid, defined by it's circuits:\\


\[
\begin{array}{cccccccccc}
n^{\circ}   & \multicolumn{9}{c}{\rm{Circuit}} \\ 
1)          & 1 & 0 & 0 & 0 & 0 & 1 & 0 & 1 & 0 \\
2)          & 0 & 1 & 0 & 0 & 1 & 0 & 0 & 1 & 0 \\
3)          & 1 & 0 & 0 & 1 & 0 & 0 & 0 & 0 & 1 \\
4)          & 0 & 1 & 1 & 0 & 0 & 0 & 0 & 0 & 1 \\
5)          & 0 & 0 & 1 & 0 & 0 & 1 & 1 & 0 & 0 \\
6)          & 0 & 0 & 0 & 1 & 1 & 0 & 1 & 0 & 0 \\
7)          & 0 & 1 & 0 & 1 & 0 & 1 & 0 & 0 & 0 \\
8)          & 1 & 0 & 1 & 0 & 1 & 0 & 0 & 0 & 0 \\
\end{array} 
\]




A possible representation should have the form:

\[
\left(
\begin{array}{ccccccccc}
1 & 0 & 0 & 1 & 1 & 1 & 1 & 1 & 0 \\
0 & 1 & 0 & x_{2,4} & 0 & x_{2,6} & x_{2,7} & x_{2,8} & 1 \\
0 & 0 & 1 & x_{3,4} & x_{3,5} & x_{3,6} & x_{3,7} & x_{3,8} & x_{3,9} \\
\end{array}
\right)
\]
The $0$ and the $1$ on the second line, depends respectively by circuits 8 and
4. The conditions expressed by the other circuits are:
\[
\left\{
\begin{array}{lr}
x_{2,6}x_{3,8}-x_{3,6}x_{2,8}=0 & (1) \\
x_{3,5}-x_{3,8}=0 & (2) \\
x_{2,4}x_{3,9}-x_{3,4}=0 & (3) \\
x_{2,6}-x_{2,7}=0 & (5)\\
x_{2,4}x_{3,5}-x_{2,4}x_{3,7}-x_{3,5}x_{2,6}+x_{3,4}x_{2,6}=0 & (6) \\
x_{3,4}-x_{3,6}=0 & (7) \\
\end{array}
\right.
\]
From this we obtain: 
\begin{equation}
\label{pappus}
x_{2,4}x_{3,5}-x_{2,4}x_{3,7}-x_{3,5}x_{2,7}+x_{3,4}x_{2,7}=0. 
\end{equation}

Now, consider the base $000000111$, this leads to 
$x_{3,9}x_{2,8}-x_{3,8}-x_{2,7}x_{3,9}+x_{3,7}\neq 0$.
By making the substitutions:
$x_{2,7}\rt x_{2,6},\ x_{3,9}\rt x_{3,4}/x_{2,4},\ 
x_{3,8}\rt x_{3,5},\ x_{2,8}\rt x_{2,6}x_{3,5}/x_{3,4}$
the preceding expression becomes: 
$-(x_{2,4}x_{3,5}-x_{2,4}x_{3,7}-x_{3,5}x_{2,7}+x_{3,4}x_{2,7})/x_{2,4}$
and this is impossible. \\

Between matroids of order $9$ and rank $3$, there are other three 
non-representable one's.
These may be obtained by adjoining other circuits to the non-Pappus matroid. 
For example $001100010$ and $110000100$. Clearly these circuits does not change 
the contradiction just showed. 
\end{example}

In the following example we show how the algorithm shows that the 
given matroid is representable only over a fields of characteristic $p$. \\

\begin{example}
\label{ex:fano}

Consider the Fano matroid:


\[
\begin{array}{cccccccc}
n^{\circ}   & \multicolumn{7}{c}{\rm{Circuit}} \\ 
1)          & 0 & 0 & 1 & 0 & 1 & 1 & 0 \\
2)          & 0 & 1 & 0 & 0 & 1 & 0 & 1 \\
3)          & 1 & 0 & 0 & 0 & 0 & 1 & 1 \\
4)          & 0 & 0 & 1 & 1 & 0 & 0 & 1 \\
5)          & 0 & 1 & 0 & 1 & 0 & 1 & 0 \\
6)          & 1 & 0 & 0 & 1 & 1 & 0 & 0 \\
7)          & 1 & 1 & 1 & 0 & 0 & 0 & 0 \\
\end{array}
\]




A possible representation should have the form:

\[
\left(
\begin{array}{ccccccccc}
1 & 0 & 1       & 0 & 1 & 0 & 1 \\
0 & 1 & x_{2,3} & 0 & 0 & 1 & x_{2,7} \\
0 & 0 & 0       & 1 & x_{3,5} & x_{3,6} & x_{3,7} \\
\end{array}
\right)
\]

The remaining conditions are:

\[
\left\{
\begin{array}{lr}
x_{2,3}x_{3,6}+x_{3,5}=0 & (1) \\
x_{3,5}-x_{3,7}=0 & (2) \\
x_{3,6}x_{2,7}-x_{3,7}=0 & (3) \\
x_{2,3}-x_{2,7}=0 & (4)\\
\end{array}
\right.
\]

This implies that $x_{2,3}x_{3,6}+x_{2,3}x_{3,6}=0$ and since the 
variables must be different from $0$, this means that the characteristic
of the field must be $2$. This leads to the well known representation of the Fano matroid 
over $\z_2$.

\[
\left(
\begin{array}{ccccccccc}
1 & 0 & 1 & 0 & 1 & 0 & 1 \\
0 & 1 & 1 & 0 & 0 & 1 & 1 \\
0 & 0 & 0 & 1 & 1 & 1 & 1 \\
\end{array}
\right)
\]

\end{example}
\label{ex:gf4}

A similar situations happens in the following example:

\begin{example}
Consider the following matroid of order $9$ and rank $4$


\[
\begin{array}{cccccccccc}
n^{\circ}   & \multicolumn{9}{c}{\rm{Circuit}} \\ 
1)           & 0 & 0 & 0 & 0 & 1 & 1 & 1 & 1 & 0 \\
2)           & 0 & 0 & 0 & 1 & 0 & 1 & 1 & 1 & 0 \\
3)           & 1 & 0 & 0 & 0 & 0 & 1 & 1 & 0 & 1 \\
4)           & 0 & 0 & 0 & 1 & 1 & 0 & 1 & 1 & 0 \\
5)           & 0 & 1 & 0 & 0 & 1 & 0 & 1 & 0 & 1 \\
6)           & 0 & 0 & 1 & 1 & 0 & 0 & 1 & 0 & 1 \\
7)           & 0 & 0 & 0 & 1 & 1 & 1 & 0 & 1 & 0 \\
8)           & 0 & 0 & 1 & 0 & 1 & 1 & 0 & 0 & 1 \\
9)           & 0 & 1 & 0 & 1 & 0 & 1 & 0 & 0 & 1 \\
10)          & 1 & 0 & 0 & 1 & 1 & 0 & 0 & 0 & 1 \\
11)          & 0 & 0 & 0 & 1 & 1 & 1 & 1 & 0 & 0 \\
12)          & 0 & 1 & 1 & 0 & 0 & 0 & 0 & 1 & 0 \\
13)          & 1 & 0 & 1 & 0 & 0 & 0 & 0 & 1 & 0 \\
14)          & 1 & 1 & 0 & 0 & 0 & 0 & 0 & 1 & 0 \\
15)          & 1 & 1 & 1 & 0 & 0 & 0 & 0 & 0 & 0 \\
\end{array}
\]




A possible representation should have the form:

\[
\left(
\begin{array}{ccccccccc}
1       & 0               & 0       & 1 & 0 & 0 & 1              & 1       & 0 \\
x_{2,1} & 0               & 1       & 0 & 1 & 0 & x_{2,1}        & x_{2,8} & 0 \\
0       & x_{2,1}x_{3,3}  & x_{3,3} & 0 & 0 & 1 & x_{2,1}x_{3,3} & x_{3,8} & 0 \\
x_{4,1} & -x_{4,3}x_{3,8} & x_{4,3} & 0 & 0 & 0 & 0              & 0       & 1 \\
\end{array}
\right)
\]

The remaining conditions are:

\[
\left\{
\begin{array}{lr}
x_{3,3}x_{2,8}-x_{3,8}+x_{2,1}x_{3,3}=0                      & (12) \\
x_{4,3}x_{3,8}+x_{4,1}x_{3,3}=0                              & (13) \\
x_{4,3}x_{2,8}-x_{2,1}x_{4,3}+x_{4,1}=0                      & (13)\\
x_{3,3}x_{2,8}-x_{3,8}-x_{2,1}x_{3,3}=0                      & (13) \\
x_{4,3}x_{3,8}-x_{4,1}x_{3,3}=0                              & (14) \\
x_{4,3}x_{2,8}x_{3,8}+x_{4,1}x_{3,8}-x_{2,1}x_{3,3}x_{4,1}=0 & (14) \\
x_{3,3}x_{2,8}+x_{3,8}-x_{2,1}x_{3,3}=0                      & (14) \\
x_{4,3}x_{2,8}+x_{2,1}x_{4,3}-x_{4,1}=0                      & (15)\\
x_{4,3}x_{2,8}+x_{2,1}x_{4,3}+x_{4,1}=0                      & (15)\\
\end{array}
\right.
\]

These equations may be satisfied only on a field of characteristic $2$, but in GF(2) there is no
solution. The first solution is in GF(4) and is given by:

\[
\left(
\begin{array}{ccccccccc}
1     & 0   & 0    & 1 & 0 & 0 & 1  & 1    & 0 \\
\e    & 0   & 1    & 0 & 1 & 0 & \e & 1+\e & 0 \\
0     & \e  & 1    & 0 & 0 & 1 & \e & 1    & 0 \\
1+\e  & \e  & 1+\e & 0 & 0 & 0 & 0  & 0    & 1 \\
\end{array}
\right)
\]

where $\e^2+\e+1=0$.

\end{example}

The equations obtained from the basis, may not admit solutions in any field, as is shown in the
following example.

\begin{example}
\label{ex:circuits}

Consider the following matroid of order $9$ and rank $4$


\[
\begin{array}{cccccccccc}
n^{\circ}   & \multicolumn{9}{c}{\rm{Circuit}} \\ 
1)           & 1 & 0 & 0 & 0 & 0 & 1 & 1 & 1 & 0 \\
2)           & 0 & 1 & 0 & 0 & 1 & 1 & 0 & 1 & 0 \\
3)           & 0 & 0 & 1 & 1 & 0 & 1 & 0 & 1 & 0 \\
4)           & 1 & 0 & 1 & 0 & 1 & 0 & 0 & 1 & 0 \\
5)           & 1 & 1 & 0 & 0 & 1 & 0 & 0 & 0 & 1 \\
6)           & 1 & 0 & 1 & 1 & 0 & 0 & 0 & 0 & 1 \\
7)           & 1 & 1 & 0 & 1 & 0 & 0 & 0 & 1 & 0 \\
8)           & 0 & 0 & 1 & 0 & 1 & 1 & 1 & 0 & 0 \\
9)           & 0 & 1 & 0 & 1 & 0 & 1 & 1 & 0 & 0 \\
10)          & 1 & 1 & 1 & 0 & 0 & 1 & 0 & 0 & 0 \\
\end{array}
\]




A possible representation should have the form:

\[
\left(
\begin{array}{ccccccccc}
1 & 0 & 1       & 1       & 0       & 0 & 1       & 0 & 1       \\
0 & 0 & x_{2,3} & 0       & 1       & 1 & x_{2,7} & 0 & x_{2,9} \\
0 & 1 & x_{3,3} & x_{3,4} & x_{3,5} & 0 & 0       & 0 & x_{3,9} \\
0 & 0 & 0       & x_{4,4} & x_{4,5} & 0 & x_{4,7} & 1 & x_{4,9} \\
\end{array}
\right)
\]

The remaining conditions are:

\[
\left\{
\begin{array}{lr}
x_{3,3}-x_{3,4}=0 & (3) \\
x_{2,3}x_{3,5}-x_{3,3}=0 & (4) \\
x_{4,9}-x_{4,5}x_{2,9}=0 & (5) \\
x_{2,3}x_{3,4}x_{4,9}-x_{2,3}x_{4,4}x_{3,9}+x_{3,3}x_{4,4}x_{2,9}=0 & (6) \\
x_{3,5}x_{4,7}+x_{3,3}x_{4,5}=0 & (8) \\
x_{4,4}-x_{4,7}=0 & (9) \\

\end{array}
\right.
\]

After substitutions, one obtain the equation
$x_{2,3}x_{4,4}x_{3,9}=0$, which is impossible.
\end{example}


\section{Results}

We have adopted algorithm \ref{alg:main} for matroids of order $8$ and $9$.
Excluded matroids are those which are non-simple or which have an element 
contained in only one circuit of cardinality less than or equal to the rank.
Matroids representable only over a finite characteristic field are listed 
in the last column.\\

The following proposition allow us to restrict our attention to matroids of 
$\rm{rank} \leq \rm{order}/2$. 

\begin{propo}
\label{pro:dual}
A matroid $M$ is representable over a field $k$ if and only if its dual is.
\end{propo}
\begin{proof}
See \cite{oxl}.
\end{proof}

The next proposition allow us to restrict our attention to matroids which do not
contain circuits of size $\leq 2$.

\begin{propo}
\label{pro:simp}
A matroid $M$ is representable over a field $k$ if and only if its
simplification $\tilde{M}$ is.
\end{propo}
\begin{proof}
See \cite{oxl}.
\end{proof}

In the following table we give a list of the results obtained by algorithm \ref{alg:main}.
In the first three columns we specify which matroids were considered and how
many of them were analyzed (third columns). In the fourth column there is the
number of non representable matroids that the algorithm has found.
In the fifth column there is the number of matroids that may be representable
only over a finite characteristic field. For these matroids we do not know if
they are representable (since our algorithm test only a sufficient condition 
for non representability) and in the affirmative case, over which field.


\begin{center}
\begin{tabular}{ccccc}
\hline
order & rank & matroids & non-rep. & finite characteristic \\
\hline
8     & 3    & 18                & 0          & 1       \\
8     & 4    & 416               & 44         & 11      \\
9     & 3    & 149               & 4          & 5       \\
9     & 4    & 179107            & 23860      & 1254    \\
10    & 3    & 2951              & 137        & 48      \\
\hline
\end{tabular}

\end{center}

Sources of the programs described in this paper are located at:
\begin{center}
{\em http://www.matapp.unimib.it/}$\sim${\em lunelli/matroidi}
\end{center}

\bibliographystyle{plain}

\end{document}